\documentclass[a4paper,12pt]{elsarticle}
\usepackage{amssymb}
\usepackage{epsfig}
\usepackage{amsfonts}
\usepackage{amsmath}
\usepackage{euscript}
\usepackage{amscd}
\usepackage{amsthm}
\DeclareMathAlphabet{\mathpzc}{OT1}{pzc}{m}{it}
\newtheorem{thm}{Theorem}[section]
\newtheorem{lem}[thm]{Lemma}
\newtheorem{prop}[thm]{Proposition}
\newtheorem{prob}[thm]{Question}

\newdefinition{defn}[thm]{Definition}
\newdefinition{ex}[thm]{Example}
\newdefinition{rem}[thm]{Remark}
\newdefinition{note}{Note}
\newdefinition{q}{Question}

\newcommand{\case}{\underline{Case}: \ }

\newcommand{\comment}[1]{}

\newcommand\m {\mathfrak{m}}
\newcommand\p {\mathfrak{p}}

\newcommand\z {\mathfrak{z}}

\begin{document}
\begin{frontmatter}
\title{Planes of the form $b(X,Y)Z^n-a(X,Y)$ over a DVR}

\author{Prosenjit Das}
\ead{prosenjit.das@gmail.com}
\author{Amartya K. Dutta}
\ead{amartya@isical.ac.in}
\address{Stat-Math Unit, Indian Statistical Institute, \\
203 B.T. Road, Kolkata 700 108, India}

\begin{abstract}
In this paper we extend an epimorphism theorem of D. Wright to
the case of discrete valuation rings. We will show that if $(R,t)$ is a
discrete valuation ring, $n \ge 2$ is an integer not divisible by the
characteristic of the residue field $R/tR$, and $g \in R[X,Y,Z]$ is a
polynomial of the form $g = b(X,Y)Z^n-a(X,Y)$ such that
$R[X,Y,Z]/(g)$ is a polynomial algebra in two variables, then $g$ and
$Z$ form a pair of variables in $R[X,Y,Z]$. We will also show that 
the result holds over any Noetherian domain containing $\mathbb{Q}$.

\end{abstract}
\end{frontmatter}

\section{Introduction} \label{sec_intro}
For a commutative ring $R$ with unity, 
let $R^{[n]}$ denote the polynomial ring in $n$ variables. 
An important question in affine algebraic geometry 
is the following epimorphism problem:

\begin{q} \label{q_1}  
Let $K$ be a field of characteristic $0$. 
Let $g \in K[X,Y,Z] (= K^{[3]})$ be such that $K[X,Y,Z]/(g) = K^{[2]}$. 
Is then $K[X,Y,Z]=K[g]^{[2]}$?
\end{q}

While the problem is open in general, a few special cases have been 
investigated by Sathaye, Russell and Wright in 
\cite{AS_OLP}, \cite{R_BIR}, \cite{DW_CNCL} and \cite{RS_FIND}; 
in some of these cases, Question 1 has an affirmative answer
even when $K$ is a field of positive characteristic.
In particular, they considered polynomials of the form 
$b(X,Y)Z^n -a(X,Y)$ and obtained affirmative answers when
\begin{enumerate}
\item [(1)] $n=1$, $K$ a field of characteristic 0 
(A. Sathaye, \cite{AS_OLP}).
\item [(2)] $n=1$, $K$ a field of any characteristic 
(P. Russell, \cite{R_BIR}).
\item [(3)] $n \ge 2$ and $K$ an algebraically closed field 
of characteristic $p \ge0$ with $p \nmid n$ (D. Wright, \cite{DW_CNCL}).
\end{enumerate}
In this paper we shall first show (see Theorem \ref{Our_Th0a})
that the above result (3) of D. Wright holds even when
$K$ is not necessarily algebraically closed.

We now consider the corresponding question over a discrete valuation ring
(to be abbreviated henceforth as DVR).

\begin{q} \label{q_2}
Let $(R, t)$ be a DVR containing $\mathbb{Q}$ and $g \in R[X,Y,Z] 
(= R^{[3]})$ be such that $R[X,Y,Z]/(g) = R^{[2]}$. 
Is then $R[X,Y,Z] = R[g]^{[2]}$?
\end{q}

As shown by Bhatwadekar-Dutta in (\cite{BD_AFNFIB}, section 4), 
this problem is closely related to the problem of $\mathbb{A}^2$-fibration 
over a regular two-dimensional
affine spot over a field of characteristic zero. 
Hence, one could explore Question \ref{q_2} at least for polynomials 
like $g = b(X,Y)Z^n-a(X,Y)$ for which the corresponding 
Question \ref{q_1} has been settled. 
For such polynomials, in view of the corresponding results over fields,
one could extend the investigation of Question 2 even to 
the positive characteristic case.

The first investigation in this direction was made by Bhatwadekar-Dutta 
in \cite{BD_DVR}. They showed (\cite{BD_DVR}, Theorem 3.5) that 
Question \ref{q_2} has an affirmative answer (in any characteristic) 
when $g=b(X,Y)Z -a(X,Y)$ with $t \nmid b(X,Y)$, thereby partially generalizing 
A. Sathaye's theorem on linear planes over a field (\cite{AS_OLP}).

The main aim of this paper is to show that Question \ref{q_2} has an
affirmative answer for polynomials of the form $g = b(X,Y)Z^n - a(X,Y)$,
where $n \ge 2$ is an integer not divisible by the characteristic of $R/tR$,
thereby obtaining a generalization of D. Wright's theorem (\cite{DW_CNCL})
quoted in section \ref{sec_prelim} (Theorem \ref{DR_Th1}).
More precisely, we will prove the following (see Theorem \ref{Our_Th4}):

\medskip

\noindent
{\bf Theorem A.} Let $(R,t)$ be a DVR with field of fractions 
$K$ and residue field $k$. Let $g \in R[X,Y,Z] (= R^{[3]})$ 
be of the form $g= bZ^n -a$ where $a,b \in R[X,Y]$ with 
$b \ne 0$ and $n$ is an integer $\ge 2$ such that $n$ is not 
divisible by the characteristic of $k$. Suppose that 
$R[X,Y,Z]/(g) = R^{[2]}$. Then $R[X,Y,Z] = {R[g, Z]}^{[1]}$, 
$R[X,Y]=R[a]^{[1]}$ and 
$b \in R[X_0]$ where $X_0 \in R[X,Y]$ and $K[X,Y] =K[X_0, a]$.

\medskip

The proof of Bhatwadekar-Dutta's theorem on linear planes over a DVR is 
highly technical. However, in the case of planes of the form $bZ^n -a$ 
with $n \ge 2$, the proof turns out to be much simpler due to the fact 
that $g$ is a variable \textit{along with} $Z$.

Using theorems on residual variables of Bhatwadekar-Dutta (\cite{BD_RES}), 
one can also see that the result for $n \ge 2$ holds over any 
Noetherian domain containing $\mathbb{Q}$. 
We shall prove (see Theorem \ref{Our_Th2d}):

\medskip

\noindent
{\bf Theorem B.} Let $R$ be a Noetherian domain containing $\mathbb Q$. 
Let $g \in R[X,Y,Z] (= R^{[3]})$ be of the form $g=bZ^n-a$ where 
$a,b \in R[X, Y]$ and $n$ is an integer $ \ge 2$. 
Suppose that $R[X,Y,Z]/(g) = R^{[2]}$. 
Then $R[X, Y, Z] = R[g, Z]^{[1]}$ and $R[X, Y] = R[a]^{[1]}$.

\medskip

In fact Theorem \ref{Our_Th2d} will show that 
the above result also holds over any Noetherian seminormal domain 
containing a field of characteristic $p \ge 0$, if $p \nmid n$.

In section \ref{sec_prelim}, we state some results which will be used
subsequently; in section \ref{sec_n=1}, we review the case $n=1$;
in sections \ref{sec_main} and \ref{sec_DVR}, we prove our main results over a
field and DVR respectively; and in section \ref{sec_misc},  we prove
our result for rings containing a field.

\section{Preliminaries} \label{sec_prelim}
Throughout this paper all rings will be commutative with unity. 
For a ring $R$, we shall use the notation 
$A = R^{[n]}$ to mean that $A$ is isomorphic, as an $R$-algebra, 
to a polynomial ring in $n$ variables over $R$; the symbol
$R^*$ will denote the group of units of $R$. 
For a prime ideal $P$ of $R$, $k(P)$ will denote the residue field $R_P/PR_P$.
An integral domain $R$ with field of fractions $K$
is called {\it seminormal} if it satisfies the condition:
an element $a \in K$ will belong to $R$ if $a^2, a^3 \in R$.

We now state some results which will be used in our proofs. 
First we state the result of D. Wright (\cite{DW_CNCL}, pg. 95) 
which we will generalize in sections 4--6.

\begin{thm} \label{DR_Th1}
Let $k$ be an algebraically closed field of characteristic $p \ge 0$. 
Let  $g \in k[X, Y, Z](= k^{[3]})$ be of the form $bZ^n-a$ 
where $a, b \in k[X, Y]$ with $b \ne 0$ and $n$ is an integer $\ge 2$ 
not divisible by $p$. Suppose that $k[X, Y, Z]/(g) = k^{[2]}$. 
Then there exist variables $\widetilde{X}, \widetilde{Y}$ in $k[X, Y]$ 
such that $a = \widetilde{Y}$, $b \in k[\widetilde{X}]$ and 
$k[X, Y, Z] = k[\widetilde{X}, g, Z]$.
\end{thm}

We now state a version of the Automorphism Theorem
of Jung and van der Kulk (\cite{Jung} and \cite{van-der-Kulk}) as presented in 
(\cite{DW_CNCL}, Appendix, Theorems 2 and 3).

\begin{thm} \label{Aut}
 Let $k$ be a field and $A = k[U,V](=k^{[2]})$. 
 Let $GA_2(k)$ denote the group of $k$-automorphisms of $A$, 
 $Af_2(k)$ the subgroup of $GA_2(k)$ defined by 
 $Af_2(k) = \{ (U,V)\mapsto (\alpha_1 U + \beta_1 V + 
 \gamma_1, \alpha_2 U + \beta_2 V + \gamma_2) | 
 \ \alpha_i, \beta_i, \gamma_i \in k \ \mbox{and} \ 
 \alpha_1\beta_2 - \alpha_2\beta_1 \ne 0 \}$, 
 ${\mathcal{E}_2}(k)$ the subgroup of $GA_2(k)$ 
 defined by ${\mathcal{E}_2}(k) = \{ (U,V) \mapsto (\alpha U + h(V), 
 \beta V + \gamma)| \ \alpha, \beta \in k^*, 
 \gamma \in k \ \mbox{and} \ h(V) \in k[V]\}$ and 
 $Bf_2(k) = Af_2(k) \cap {\mathcal{E}_2}(k)$. 
 Then $GA_2(k) = Af_2(k) *_{Bf_2(k)} \mathcal{E}_2(k)$. 
 Moreover, if $\sigma \in GA_2(k)$ is of finite order, 
 then there exists $\tau \in GA_2(k)$ such that either 
 $\tau \sigma \tau^{-1} \in Af_2(k)$ or 
 $\tau \sigma \tau^{-1} \in\mathcal{E}_2(k)$.
\end{thm}

Now we state a result of A. Sathaye (\cite{AS_OLP}, Corollary 1) 
which we will use to prove Lemma \ref{Our_Lem2}.

 \begin{thm} \label{Sat_lem}
Let $L|_k$ be a separable field extension. Assume that there exist 
$h \in k[X,Y]$ and  $f_i \in L[X,Y]$, $1 \le i \le s$, such that
\begin{enumerate}
\item [\rm \textit{(1)}]$L[X,Y]/(f_i) = L^{[1]}$ for each $i$.
\item [\rm \textit{(2)}]$(f_i, f_j)L[X,Y] = L[X,Y]$ for $i \ne j$.
\item [\rm \textit{(3)}]$h = \underset{i=1} {\overset{s} {\prod}} 
{f_i}^{r_i}$, $r_i > 0$.
\end{enumerate}
Then there exist $f \in k[X,Y]$, $\lambda_i \in L^*$ and $\mu_i \in L$ 
such that $f_i = \lambda_i f+ \mu_i$ for each $i$, $1 \le i \le s$.
\end{thm}

We will also use the following special case of the result
 (\cite{D_SEP}, Theorem 7).

\begin{thm} \label{Dut_Sep}
Let $k$ be a field, $L$ a separable field extension of $k$, 
$A$ a UFD containing $k$ and $B$ an $A$-algebra such that 
$B \otimes_k L = (A \otimes_k L)^{[1]}$. Then $B=A^{[1]}$.
\end{thm}

We will use the following version of a cancellation theorem due to
 Abhyankar-Eakin-Heinzer (\cite{AEH_Coff}, Theorem 3.3).

\begin{thm} \label{AEH}
Let $A$ be an affine domain over a field $k$ such that $k$ is 
algebraically closed in $A$ and $tr. deg_k(A) =1$. 
Suppose that $B$ is another $k$-algebra such that $A^{[n]} = B^{[n]}$ 
for some $n \ge 1$. 
Then either $B=A$ or $B \cong A = k^{[1]}$.
\end{thm}

We now state a version of the Russell-Sathaye criterion 
(\cite{RS_FIND}, Theorem 2.3.1) for a ring to be a polynomial algebra 
over a subring (see \cite{BD_DVR}, Theorem 2.6).

\begin{thm} \label{RS_Poly}
Let $R \subset A$ be integral domains with $A$ being finitely generated 
over $R$. Suppose that there exist primes $p_1, p_2, \dots , p_n$ in 
$R$ such that for each $i$, $1 \le i \le n$,   
\begin{enumerate}
\item[(1)] $p_i$ remains prime in $A$, 
\item[(2)] $p_iA \cap R = p_iR$,
\item[(3)] $A[\frac{1}{p_1 p_2 \dots p_n}] = 
R[\frac{1}{p_1 p_2 \dots p_n}]^{[1]}$ and
\item[(4)] $R/p_iR$ is algebraically closed 
in $A/p_iA$.
\end{enumerate}
Then $A = R^{[1]}$.
\end{thm}

The following result from (\cite{BD_DVR}, 2.5) will enable 
us to apply Theorem \ref{RS_Poly}.

\begin{lem} \label{Lem_Alg_Closed}
 Let $R$ be an integral domain and
 $F \in R[X, Y](=R^{[2]})$ be such that
 $R[X, Y]/(F) = R^{[1]}$.
 Then $R[F]$ is algebraically closed in $R[X,Y]$.
\end{lem}

Finally, we state a result on residual variables which will be
our main tool to prove Theorem B.
It comes as a direct consequence of Theorem 3.1, 
Theorem 3.2 and Remark 3.4 in \cite{BD_RES}. 

\begin{thm} \label{Res_Result}
Let $R$ be a Noetherian domain such that either
$R$ contains $\mathbb{Q}$ or $R$ is seminormal,
$A$ be a polynomial algebra in $n$ variables over $R$ and $W_1, W_2, \dots , W_{n-1} \in A$. Then the following are equivalent:

1. $A = R[W_1, W_2, \dots , W_{n-1}]^{[1]}$.

2. $A \otimes_R k(P) = (R[{W_1}, {W_2}, \dots , {W_{n-1}}] \otimes_R k(P))^{[1]}$ for every prime ideal $P$ of $R$.
\end{thm}

\section{Planes of the form $bZ-a$} \label{sec_n=1}

We recall below the earlier result on linear planes over a DVR
(\cite{BD_DVR}, Theorem 3.5).

\begin{thm} \label{Dut_DVR_Th1}
 Let $(R,t)$ be a DVR and $g \in R[X,Y,Z](= R^{[3]})$ be of the form
 $g=bZ -a$ where $a,b \in R[X,Y]$ and $b \notin tR[X,Y]$.
 Suppose that $R[X,Y,Z]/(g) = R^{[2]}$. Then $R[X,Y,Z] = R[g]^{[2]}$.
\end{thm}

We now show that the result can be generalized to the case of
Dedekind domain in the following form.

\begin{thm} \label{Dut_DVR_Th2}
 Let $R$ be a Dedekind domain and $g \in R[X,Y,Z](= R^{[3]})$ be of the form
 $g= bZ -a$ where $a,b \in R[X,Y]$ and the coefficients of $b$ generate the
 unit ideal of $R$. Suppose that $B=R[X,Y,Z]/(g) = R^{[2]}$.
 Then $R[X,Y,Z] = R[g]^{[2]}$.
\end{thm}

\begin{proof}
By Theorem \ref{Dut_DVR_Th1}, $R_{\m}[X,Y,Z] = R_{\m}[g]^{[2]}$ for 
each maximal ideal $\m$ of $R$. Hence, by (\cite{BCR_Local}), 
it follows that $R[X,Y,Z]$ is $R[g]$-isomorphic to the symmetric 
algebra $Sym_{R[g]}(P)$ for some finitely generated projective 
$R[g]$-module $P$ of rank two. Thus it is enough to show that $P$ 
is a free $R[g]$-module. Since $R[g]$ is a retract of $R[X,Y,Z]$, 
it is enough to show that $P \otimes_{R[g]} R[X,Y,Z]$ is a free 
$R[X,Y,Z]$-module. Note that since $R[X,Y,Z] \cong Sym_{R[g]}(P)$, 
we have $\Omega_{R[g]}(R[X,Y,Z]) = P \otimes_{R[g]} R[X,Y,Z]$. 
Thus the proof will be complete if we show that the projective 
$R[X,Y,Z]$-module $\Omega_{R[g]}(R[X,Y,Z])$ is actually free.

Now consider the exact sequence:
\begin{eqnarray*}
\Omega_R (R[g]) \otimes_{R[g]} (R[X,Y,Z]) \underset{} {\overset{\theta} {\longrightarrow}} \Omega_R (R[X,Y,Z]) \longrightarrow \Omega_{R[g]} (R[X,Y,Z]) \longrightarrow 0.
\end{eqnarray*}
Let $g_X$, $g_Y$ and $g_Z$ denote the partial derivatives of $g$ 
with respect to $X$, $Y$ and $Z$ respectively. Now note that 
$(g_X, g_Y, g_Z)R[X,Y,Z] = R[X,Y,Z]$. Since $dim \ R = 1$, by 
Suslin's theorem (\cite{Sus_SL}, 2.6), the unimodular row 
$[g_X, g_Y, g_Z]$ can be completed to an invertible matrix. 
Since $\Omega_R(R[X,Y,Z])$ is a free $R[X,Y,Z]$-module of rank 
three with basis $dX$, $dY$ and $dZ$, and since $Im \ (\theta)$ 
is generated by $g_X dX + g_Y dY + g_Z dZ$, it now follows that 
$\Omega_{R[g]}(R[X,Y,Z]) (= \Omega_R (R[X,Y,Z])/ Im \ (\theta))$ 
is a free $R[X,Y,Z]$ module of rank two. This completes the proof.
\end{proof}

\begin{rem}
Let $(R,t)$ be a DVR containing $\mathbb{Q}$ and 
let $g = bZ-a$ where $b = tY^2$ and $a = -Y - tY(X+X^2) - t^2 X$. 
Then $R[X,Y,Z]/(g) = R^{[2]}$ (see \cite{BD_AFNFIB}, Example 4.13).
In this example, $t \mid b$; and 
it is not yet known whether $R[X,Y,Z] = R[g]^{[2]}$. 
\end{rem}

\section{Planes of the form $bZ^n-a$ over a field} \label{sec_main}

In this section we will show that Wright's arguments in (\cite{DW_CNCL})
 can be modified to show that his result (Theorem \ref{DR_Th1}) can be extended over
 any field. We first prove a few auxiliary results 
(Lemmas \ref{auto} and \ref{Our_Lem2}), then consider
 the case when the field $k$ contains all $n$th roots of unity
 (Proposition \ref{Our_Th0}) and finally show that Theorem \ref{DR_Th1}
 holds over any field (Theorem \ref{Our_Th0a}).
We first record a result on ${\mathrm Aut}_k (k^{[2]})$.

\begin{lem} \label{auto}
Let $k$ be a field of characteristic $p \ge 0$ and $\sigma$ a 
$k$-automorphism of $B = k^{[2]}$ of order $n$ such that 
$p \nmid n$. Suppose that $k$ contains all the $n^{th}$ roots 
of unity. Then there exist elements $U, V \in B$ and 
$\alpha, \beta \in k^*$ such that $B = k[U, V]$, $\sigma (U)= \alpha U$ 
and $\sigma (V) =\beta V$, where $\alpha^n = \beta^n  = 1$.
\end{lem}

\begin{proof}
By Theorem \ref{Aut}, one can choose coordinates $U'$, $V'$ of $B$ 
such that either $\sigma \in {\mathcal{E}_2}(k)$ or $\sigma \in Af_2(k)$.

\smallskip

\case $\sigma \in \mathcal{E}_2(k)$.

In this case $\sigma (U') = \alpha U' + \mu$ and 
$\sigma(V') =\beta V' + f_1(U')$, 
where $\alpha, \beta \in k^*$, $\mu \in k$ and $f_1(U') \in k[U']$. 
Since $\sigma$ is of order $n$, we have $\alpha^n = \beta^n =1$. Note that 
if $\alpha = 1$, then $U' = \sigma^n(U') = U' + n \mu$ and hence 
$\mu = 0$, as $p \nmid n$. 

Set 
$$
U: = \left \{ \begin{array}{lll}
		    U' & \mbox{if} & \alpha = 1 .\\
		    U' + \frac{\mu}{\alpha - 1} & \mbox{if} & \alpha \ne 1.
                 \end{array}
	\right.
$$
 Then 
$k[U',V'] = k[U, V']$, $\sigma(U) = \alpha U$ and 
$\sigma(V') = \beta V' + f(U)$ 
for some $f(U) \in k[U]$. We will now show that we can choose 
$g(U) \in k[U]$ such that $\sigma (V' + g(U)) = \beta (V' + g(U))$. 
Let $f(U) = \underset{i=0} {\overset{r} {\varSigma}} a_i {U}^i$. 

First we show that for any $i$, $1 \le i \le r$, if $a_i \ne 0$, 
then $\alpha^i \ne \beta$. Suppose $\beta = \alpha^i$. Now, from 
the relation $V' = \sigma^n (V')$, we get 
\[
\beta^{n-1} f(U) + \beta^{n-2}f(\alpha U) + \dots + f(\alpha^{n-1} U) = 0,
\] 
which implies that 
\[
\beta^{n-1} a_i + \beta^{n-2} \alpha^i a_i + \beta^{n-3}\alpha^{2i}a_i + 
\dots + \alpha^{(n-1)i} a_i = 0,
\] 
i.e., $n \beta^{n-1} a_i = 0$, and hence $a_i =0$ 
(as $p \nmid n$ and $\beta \ne 0$). Thus $\alpha^i \ne \beta$ if $a_i \ne 0$.

Now, for each $i$, $ 1 \le i \le r$, we define $b_i$ as follows:
$$
b_i = \left \{ \begin{array}{lll}
		    0 & \mbox{if} & a_i = 0.\\
		    a_i/(\beta - \alpha^i) & \mbox{if} & a_i \ne 0.
                 \end{array}
	\right.
$$
Let $g(U) = \underset{i=0} {\overset{r} {\varSigma}} b_i {U}^i$ 
and set 
\[
V := V' + g(U).
\]
Then $\sigma(V)=\beta V$.
Thus $k[U',V'] = k[U,V]$, $\sigma(U) = \alpha U$ and $\sigma (V) = \beta V$.

\smallskip

\case $\sigma \in  Af_2(k)$.

In this case $\sigma (U')=\alpha_1 U' + \beta_1 V' + \gamma_1$ and 
$\sigma(V')= \alpha_2 U' + \beta_2 V' + \gamma_2$ 
for some $\alpha_i, \beta_i, \gamma_i \in k$ ($i=1,2$) 
with $\alpha_1 \beta_2 \ne \beta_1 \alpha_2$. 
Choose $\lambda \in \bar{k}$ such that 
$(\alpha_1 - \lambda) (\beta_2 - \lambda) - \alpha_2 \beta_1 =0$. 
Then $\lambda$ is an eigen value of the linear transformation 
$(X,Y) \mapsto (\alpha_1 X + \alpha_2 Y, \beta_1 X + \beta_2 Y)$ of 
$\bar{k}^2$. 
Let $(\nu_1, \nu_2) \in \bar {k}^2$, not both zero, be an eigen vector 
corresponding to the eigen value $\lambda$. Then we have
\begin{eqnarray*}
\alpha_1 \nu_1 + \alpha_2 \nu_2 &=& \lambda \nu_1 \\ \beta_1 \nu_1 + \beta_2 
\nu_2 &=& \lambda \nu_2.
\end{eqnarray*}
Therefore,
$\sigma(\nu_1 U' + \nu_2 V') = \lambda (\nu_1 U' + \nu_2 V') + \mu$
where $\mu=\nu_1 \gamma_1 + \nu_2 \gamma_2$.
Since $\sigma$ is of order $n$, we have 
$\lambda^n = 1$ and hence $\lambda \in k^*$. 
Thus we may choose $\nu_1, \nu_2 \in k$. 
Therefore, setting $U := \nu_1 U' + \nu_2 V'$, 
we have $\sigma(U) = \lambda U + \mu$ and hence 
$\sigma(V') = \kappa V' + h(U)$ for some 
$\kappa \in k^*$ and $h(U) \in k[U]$. Now, by taking $U$ and $V'$ 
to be the coordinates for $B$, the problem reduces to the 
previous case: $\sigma \in \mathcal{E}_2(k)$.

\smallskip

Thus in both the cases we get $U, V \in B$ and $\alpha, \beta \in k^*$ 
such that $B = k[U, V]$, $\sigma (U) = \alpha U$ and $\sigma(V) = \beta V$. 
This completes the proof.
\end{proof}

We now record a consequence of Sathaye's result (Theorem \ref{Sat_lem}).

\begin{lem} \label{Our_Lem2}
 Let $k$ be a field, $B = k^{[2]}$ and $b \in B \backslash k$. 
 Suppose that there exist a separable algebraic extension $E|_k$ 
 and an element $X' \in B \otimes_k E$ such that 
 $B \otimes_k E = E[X']^{[1]}$ and $b \in E[X']$. 
 Then there exists $X \in B$ such that $b \in k[X]$, 
 $B = k[X]^{[1]}$ and $E[X']=E[X]$.
\end{lem}

\begin{proof}
Without loss of generality, we assume $E|_k$ to be a finite 
Galois extension. 
Let $B=k[X_1, Y_1]$. Then $B \otimes_k E= E[X_1, Y_1]=E[X']^{[1]}$.
Let $X'=\phi(X_1,Y_1)$. Interchanging $X_1$ and $Y_1$ if necessary,
we may assume that the $X_1$-degree of $\phi(X_1,Y_1)$ is positive.
Hence the leading coefficient of $X_1$ in $\phi(X_1,Y_1)$ is a non-zero
element $\lambda \in E$ (\cite{Abh_LecNote_Expan}, Proposition 11.12, pg. 85). 
Let $X''=X'/\lambda$.

Let $G = \{ \sigma_i \ | \ i=1,2, \dots , m \}$ be the group of 
$k$-automorphisms of $E|_k$. We extend each $\sigma \in G$ 
to a $B$-automorphism of $B \otimes_k E$. 
Let $\bar{k}$ be an algebraic closure of $k$ containing $E$ and 
$b=\underset{i=1} {\overset{s} {\prod}} (\lambda_i X'' + \mu_i)^{n_i}$ 
be the prime decomposition of $b$ in $\bar{k}[X'']$, where 
$\lambda_i \in {\bar{k}}^{^{^{*}}}$, $\mu_i \in \bar{k}$ and 
$n_i \in \mathbb{N}$, 
$1 \le i \le s$. Since $\sigma(b) = b$ for each $\sigma \in G$, 
$b = \underset{i=1} {\overset{s} {\prod}} (\sigma(\lambda_i) 
\sigma(X'') + \sigma(\mu_i))^{n_i}$ is also a prime decomposition of $b$ in 
$\bar{k}[X'']$. This shows that for each $\sigma \in G, \ 
\exists \alpha \in {\bar{k}}^{^{^{*}}}$ and $\beta \in \bar{k}$ such that 
$\sigma(X'') = \alpha X'' + \beta$. Since $X''$ and $\sigma(X'')$
are both monic in $X_1$, it follows that $\alpha =1$.

Since $X''$ is a variable of $B \otimes_k E$, we have 
$(B \otimes_k E)/(\sigma(X'')) = E^{[1]}$ for each $\sigma \in G$. 
It is also easy to see that if $\sigma_i(X'') \ne \sigma_j(X'')$ 
for $\sigma_i, \sigma_j \in G$,
then $\sigma_i(X'')$ and $\sigma_j(X'')$ are comaximal in 
$B \otimes_k \bar{k}$ and hence comaximal in $B \otimes_k E$. 
Let $f_1, \dots , f_t$ be the distinct elements of the set 
$\{ \sigma(X'')| \sigma \in G\}$.  Then, for each $i$, $1 \le i \le t$, 
there exists $m_i \in \mathbb{N}$ such that 
$\underset{ \sigma \in G} {\overset{} 
{\prod}} \sigma(X'') = \underset{i=1} {\overset{t} {\prod}} 
f_i^{m_i} \in B$, $(B \otimes_k E)/(f_i) = E^{[1]}$, and for 
$i \ne j$, $f_i$ and $f_j$ are comaximal in $B \otimes_k L$. 
Since $B = k^{[2]}$, applying Theorem \ref{Sat_lem}, we get that for 
each $\sigma \in G$ there exist $\lambda \in E^*$ and $\mu \in E$ 
such that $\lambda \sigma(X'')+ \mu \in B$. Fix $\sigma \in G$ and 
let $X = \lambda \sigma(X'')+ \mu \in B$. Then  $E[X'']=E[\sigma(X'')]=E[X]$ 
and $b \in E[X] \cap B$. Since $B = k^{[2]}$ and $X \in B$, 
we have $E[X] \cap B = k[X]$. Hence 
$b \in k[X] \subset B$. Now since $B = k^{[2]}$ and 
$B \otimes_k E = E[X'']^{[1]}= E[X]^{[1]}$, by 
Theorem \ref{Dut_Sep}, we see that $B = k[X]^{[1]}$. By 
construction, $E[X]=E[X'']=E[X']$.
\end{proof}

For convenience, we state below a result
which follows from a lemma of A. Sathaye (\cite{AS_OLP}, Lemma 1).

\begin{lem} \label{Sat_lem1}
Let $k$ be a field and suppose $X'$ is a variable in 
$k[X_1, X_2, \dots, X_n] (= k^{[n]})$ which is comaximal with $X_1$. 
Then $X' = \alpha X_1 + \beta$ with $\alpha, \beta \in k$, $\alpha \ne 0$.
\end{lem}

\begin{prop}\label{Our_Th0}
Let $k$ be a field of characteristic $p \ge 0$ containing the $n^{th}$ 
roots of unity and $g \in k[X,Y,Z](= k^{[3]})$ be of the form 
$bZ^n -a$ where $a, b \in k[X,Y]$ with $b \ne 0$ and $n$ is an 
integer $\ge 2$ not divisible by $p$. 
Suppose that $B := k[X,Y,Z]/(g) = k^{[2]}$. 
Then there exist variables $U,V$ in $B$ such that $V$ is the 
image of $Z$ in $B$, $b \in k[U]$ and $k[X,Y] = k[U,a] = k^{[2]}$.
\end{prop}

\begin{proof}
Let $ \sigma$ be the $k$-automorphism of $B$ induced by the 
$k$-automorphism $\widetilde{\sigma}$ of $k[X,Y,Z]$ defined by 
$\widetilde{\sigma}((X,Y,Z)) = (X,Y, \omega Z)$ where $\omega$ is a 
primitive $n^{th}$ root of unity. Obviously, $\sigma$ has order $n$.

Since $B = k^{[2]}$, by Lemma \ref{auto}, there exist elements 
$U', V' \in B$ and $\alpha, \beta \in k^*$ such that $B = k[U', V']$, 
$\sigma(U')=\alpha U'$ and $\sigma(V')=\beta V'$, where 
$\alpha^n = \beta^n  = 1$. Let $\z$ be the image of $Z$ in $B$ 
and $A = k[X,Y][a/b]$. Then $\z^n = a/b$ and 
$B = A[\z] = k[X,Y][\z] = A \oplus \z A \oplus \z^2 A \oplus \dots \oplus \z^{n-1}A$ 
so that, for any $x \in B$, $\z \mid (x-\sigma(x))$. 
Thus $\z \mid (1-\alpha)U'$ 
and $\z \mid (1-\beta)V'$. But since $U'$ and $V'$ can not have 
common (non-unit) factor and $\z \notin k^*$, we have either 
$\alpha = 1$ or $\beta = 1$. Interchanging $U'$ and $V'$ if 
necessary, we assume that $\alpha = 1$. Then the ring of 
invariants of $\sigma$ is $A = k[X,Y][a/b]=k[U', a/b] (= k^{[2]})$. 
Note that $V'$ is a unit multiple of $\z$. 
Thus $B=k[U', \z]$. Set $V:=\z$.

Now we show that we can choose $U$ from $k[X,Y]$ 
such that $B=k[U,V]$, 
$b \in k[U]$ and $k[X,Y] = k[U,a]$. 
If $b \in k^*$, 
then $k[X,Y]=k[X,Y][a/b]=k[U',a/b]$, so that, in this case, 
we may set $U:=U'$.  We now consider the case 
 $b \notin k^*$. Let $p_1, p_2, \dots , p_m$ be the distinct 
irreducible factors of $b$ in $A (=k^{[2]})$, and set
${\p}_i: = k[X,Y] \cap p_i A$. Note that for each $i = 1,2, \dots, m$, 
both $b$ and $a (=b.a/b) \in k[X,Y] \cap bA \subseteq {\p}_i$. 
This shows that $(bZ^n -a)k[X,Y,Z] \subsetneqq {\p}_i[Z]$ which implies 
$ht \ {\p}_i > 1$. Thus each ${\p}_i$ is a maximal ideal of $k[X,Y]$. 
Let $\bar{k}$ denote an algebraic closure of $k$, $L_i$ be a
subfield of $\bar{k}$ isomorphic to $k[X,Y]/\p_i$ and let $L$ be  
the subfield of $\bar{k}$ generated by the fields $L_1, L_2, \dots, L_m$. 
 Then $L_i$ is an algebraic extension of $k$ 
and  $A/\p_i A = (k[X,Y]/\p_i) [\zeta_i] = L_i[\zeta_i]$ where 
$\zeta_i$ is the image of $a/b$ in $A/\p_i A$. Since 
$\p_iA \subseteq p_i A$, it follows that 
$\zeta_i$ is transcendental over $L_i$ and 
$\p_iA$ is a prime ideal of $A$. As $ht \ p_iA =1$ and  
$\p_i A \ne 0$, we have $p_i A = \p_i A$. This shows that $p_i$ are pairwise 
comaximal in $A$ and hence in $B$.

Let $g(\zeta_i)$ be the image of $U'$ in $A/p_i A = L_i[\zeta_i]$. 
Then $U' -g(a/b)$ is divisible by $p_i$ in $A \otimes_k L_i$. But 
$U' - g(a/b) = U'-g(V^n)$ is a variable in both $A\otimes_k L_i$ and 
$B \otimes_k L_i$. Hence $U' - g(a/b)$ is a constant multiple of $p_i$. 
Thus $A \otimes_k L_i=L_i[p_i, a/b]$, 
$B \otimes_k L_i = L_i[p_i, V]$, and for 
$i \ne j$, $(p_i, p_j)B \otimes_k L = B \otimes_k L$. 
Set $U:=p_1$.
Using Lemma \ref{Sat_lem1}, we have 
$p_i = \lambda_i U + \mu_i$ for $\lambda_i \in L^*$ and $\mu_i \in L$. 
So, we have $b \in L[U]$. This shows that $U$ is integral 
over $L[X,Y]$ and hence over $k[X,Y]$.
As $U \in k[X,Y][a/b]$ and $k[X,Y]$ is a normal domain, 
we have $U \in k[X,Y]$. 
Since $L|_k$ is faithfully flat, it follows that $B = k[U,V]$ with 
$U \in k[X,Y], V = \z$ and $b \in k[U]$. Now, the argument in 
(\cite{DW_CNCL}, pg. 98) shows that $k[X,Y] = k[U,a]$.
\end{proof}

\begin{thm} \label{Our_Th0a}
 Let $k$ be a field of characteristic $p \ge 0$ and 
$g \in k[X,Y,Z]$ be of the form $bZ^n -a$ where 
$a, b \in k[X,Y]$ with $b \ne 0$ and $n$ is an integer 
$\ge 2$ not divisible by $p$. Suppose that 
$B := k[X,Y,Z]/(g) = k^{[2]}$ and identify $k[X,Y]$
with its image in $B$. Then there exist variables 
$U, V$ in $B$ such that $V$ is the image of $Z$ in $B$, $U \in k[X,Y]$,
$b \in k[U]$, $k[X,Y] = k[U,a]$ and $k[X,Y,Z] = k[U,g,Z]$.
\end{thm}

\begin{proof}
Let $E$ be the field obtained by adjoining 
all the $n^{th}$ roots of unity to $k$.
Since $p \nmid n$, $E$ is a Galois extension over $k$.
By Proposition \ref{Our_Th0}, we get variables $U'$ and 
$V'$ of $B \otimes_k E$ ($=k[X,Y,Z]/(g) = E^{[2]}$) 
 such that $V'$ is the image of $Z$, 
 $b \in E[U']$ and $E[X,Y] = E[U',a]$. 
As $E|_k$ is separable, we have $k[X,Y] = {k[a]}^{[1]}$ by
Theorem \ref{Dut_Sep}. 
If $b \in k[X,Y] \backslash k$, then, by Lemma \ref{Our_Lem2},
we get $U \in k[X,Y]$ such that $k[X,Y] = k[U]^{[1]}$, $b \in k[U]$ and 
$E[U] = E[U']$. Since $E|_k$ is faithfully flat,
 $E[U',a] = E[U, a]$ and $k[U,a] \subseteq k[X,Y]$, we have 
$k[U, a]=k[X,Y]$. 
If $b \in k$, then we choose $U$ to be any complementary variable
of $a$ in  $k[X,Y]$. 

From the relation  $k[U, a]=k[X,Y]$, we have
\[
k[X,Y,Z] = k[U,a,Z] = k[U, bZ^n -a, Z] = k[U,g,Z].
\]
The relation $k[X,Y,Z]=k[U,g,Z]$ shows that $B$ is generated
by the images of $U$ and $Z$. This completes the proof.
\end{proof}

\begin{rem}
Theorem \ref{Our_Th0a} does not hold if $p \mid n$. Consider a 
field $k$ of characteristic $p > 0$ and the polynomial
$g = Z^{p^e} - Y - Y^{sp} \in k[Y,Z]$ where $p \nmid s$ and $e \ge 2$. Then 
$\frac{k[Y,Z]}{(g)} = k^{[1]}$ but $k[Y,Z] \ne k[g]^{[1]}$ 
(see \cite{Abh_LecNote_Expan}, Example 9.12, pg. 72). 
Using a result of E. Hamann (\cite{Haman_Invariance}, Theorem 2.6), 
it follows that $k[X,Y,Z] \ne k[g]^{[2]}$ although 
$\frac{k[X,Y,Z]}{(g)} = k^{[2]}$.
\end{rem}

\section{Planes of the form $bZ^n-a$ over a DVR} \label{sec_DVR}
In this section we shall prove Theorem A. 
We first record two results on factorial domains.

\begin{lem} \label{Lem_Inert}
 Let $R$ be a UFD with field of fractions $K$. Let $U \in R[X,Y]$ 
be such that $K[X,Y] = K[U]^{[1]}$. Then $K[U] \cap R[X,Y]$ is an 
inert subring of $R[X,Y]$ and $K[U] \cap R[X,Y] =R[W](= R^{[1]})$, where
$W$ is an element of $R[X,Y]$ such that $K[W]=K[U]$.

\end{lem}
\begin{proof}
 Let $D = K[U] \cap R[X,Y]$. Clearly, $D$ is an inert subring of 
$R[X,Y]$ and hence a UFD of transcendence degree one over $R$. Therefore,
by (\cite{AEH_Coff}, Theorem 4.1), 
$D = R[W] (= R^{[1]})$ for some $W \in R[X,Y]$. Clearly, $K[W]=K[U]$.
\end{proof}

\begin{lem} \label{Lem_Inert_R[a]}
Let $R$ be a UFD of characteristic $p \ge 0$ with field of fractions $K$ 
and $g \in R[X,Y,Z] (= R^{[3]})$ be of the form $g= bZ^n -a$ where 
$a,b \in R[X,Y]$ with $b \ne 0$ and $n$ is an integer $\ge 2$ 
such that $p \nmid n$. Suppose that $R[X,Y,Z]/(g) = R^{[2]}$. Then
\begin{enumerate}
\item [\rm (i)]$R[a]=K[a] \cap R[X,Y]$.
\item [\rm (ii)]$R[a]$ is an inert subring of $R[X,Y]$.
\item [\rm (iii)]$tR[X,Y] \cap R[a] = tR[a]$ for every $t \in R$.
\end{enumerate}
\end{lem}
\begin{proof}
(i) \ By Theorem \ref{Our_Th0a}, $K[X,Y] = K[a]^{[1]}$ and
by Lemma \ref{Lem_Inert}, $K[a] \cap R[X,Y] = R[W]$ for some 
$W \in R[X,Y]$ satisfying $K[a] = K[W]$. It then follows that 
$a = \lambda W + \mu$ where $\lambda, \mu \in R$. 
We claim that $\lambda \in R^*$. Suppose $\lambda \notin R^*$. 
Let $q$ be a prime factor of $\lambda$ and let
$L$ denote the algebraic closure of the field of fractions of $R/qR$. Let
$\bar{a}$ and $\bar{b}$ denote
the images of $a$ and $b$ respectively in $L[X,Y]$.
Then we would have $\bar{a}(=\mu) \in L$; in fact, as 
 $L[X,Y,Z]/(g) = L[X,Y,Z]/(\bar{b}Z^n - \bar{a}) = L^{[2]}$,
 we would have that $\bar{a}$ is a unit in $L$. Since 
$L[X,Y] \hookrightarrow L[X,Y,Z]/(\bar{b}Z^n - \bar{a}) (= L^{[2]})$,
it would follow that $\bar{b} \in L^*$. But then, as $n \ge 2$, 
$L[X,Y,Z]/(\bar{b}Z^n - \bar{a})$ would not be an integral domain,
contradicting that $L[X,Y,Z]/(\bar{b}Z^n - \bar{a}) = L^{[2]}$. 
Thus $\lambda \in R^*$ and hence $R[a]=R[W]=K[a] \cap R[X,Y]$.

(ii) and (iii) follow from (i).
\end{proof}

We now prove Theorem A.

\begin{thm} \label{Our_Th4}
 Let $(R, t)$ be a DVR with residue field $k$ and let $p (\ge 0)$ be
 the characteristic of $k$. Let $g \in R[X,Y,Z] (= R^{[3]})$ be of the form
 $g= bZ^n -a$ where $a,b \in R[X,Y]$ with $b \ne 0$ and $n$ is an
 integer $\ge 2$ such that $p \nmid n$.
 Suppose that $R[X,Y,Z]/(g) = R^{[2]}$.
 Then  $R[X,Y,Z] = R[g, Z]^{[1]}$, $R[X,Y]=R[a]^{[1]}$
 and $b \in R[X_0]$
 where $X_0 \in R[X,Y]$ and $K[X,Y]=K[X_0,a]$.
\end{thm}

\begin{proof}
Let $K$ and $k$ denote, respectively, the field of fractions and 
the residue field 
of $(R,t)$. For any $f \in R[X,Y,Z]$, let $\bar{f}$ denote 
the image of $f$ in $k[X,Y,Z]$. By hypotheses, $K[X,Y,Z]/(bZ^n -a) = K^{[2]}$
and $k[X,Y,Z]/(\bar{b}Z^n -\bar{a}) = k^{[2]}$. 
Hence, by Theorem \ref{Our_Th0a}, $K[X,Y]=K[a]^{[1]}$ 
and $K[X,Y,Z] = K[Z, bZ^n -a]^{[1]}$.

If $t \nmid b$, then, by Theorem \ref{Our_Th0a}, 
$k[X,Y,Z] = k[Z, \bar{g}]^{[1]}$ and
$k[X,Y]=k[\bar{a}]^{[1]}$.  Hence,
by Theorem \ref{RS_Poly}, we get 
$R[X,Y,Z] = R[g,Z]^{[1]}$ and $R[X,Y]=R[a]^{[1]}$.

We now consider the case $t \mid b$. Now $\bar{g}=\bar{a}$ so that 
\[
 k[X,Y,Z]/(\bar{a}) \ (=k[X,Y]/(\bar{a}))^{[1]} \ 
= k[X,Y,Z]/(\bar{g}) = k^{[2]}.
\]
Hence, by Theorem \ref{AEH}, $k[X,Y]/(\bar{a}) = k^{[1]}$.
Therefore, by Lemma \ref{Lem_Alg_Closed}, 
we see that $k[\bar{a}]$ is algebraically closed in $k[X,Y]$. 
Since $t$ is prime in both $R[a](=R^{[1]})$ and $R[X,Y]$,
and since $a$ is a generic variable of $R[X,Y]$, using Theorem \ref{RS_Poly},
we see that $R[X,Y] = R[a]^{[1]}$. By similar argument, we have 
$R[X,Y,Z] = R[g, Z]^{[1]}$.

Now, by Theorem \ref{Our_Th0a}, one can choose $U \in R[X,Y]$ such that
$K[X,Y] = K[U,a]$ and $b \in K[U]$. By Lemma \ref{Lem_Inert}, 
$K[U] \cap R[X,Y] = R[X_0]$ for some $X_0 \in R[X,Y]$ satisfying
$K[U] = K[X_0]$. Thus $b \in R[X_0]$ 
where $K[X_0, a]=K[U,a]=K[X,Y]$. Hence the result.
\end{proof}

Note that, in the case $R$ is a $\mathbb{Q}$-algebra, 
the hypothesis in Theorem \ref{Our_Th4} regarding $n$ 
($p \nmid n$) is automatically satisfied.
Thus, in particular, Theorem \ref{Our_Th4} holds when $R$ is a
DVR containing $\mathbb{Q}$. In the next section
we shall see a generalisation of this result
(Theorem \ref{Our_Th2d}).

\begin{rem}
Note that, in the notation of Theorem \ref{Our_Th4},
 $X_0$ need not be a variable in $R[X,Y]$. Consider a DVR $(R,t)$. Let
$g = bZ^n - a$ where $a = -Y$ and $b = t^2 X + t Y^2$, and let 
$X_0 = t X + Y^2$. Then $R[X,Y,Z]/(g) = R^{[2]}$, 
$b \in R[X_0]$, $K[X,Y] = K[X_0, Y]$ but $R[X,Y] \ne R[X_0]^{[1]}$.
\end{rem}

The following example shows that, without the hypothesis $p \nmid n$, 
Theorem \ref{Our_Th4} need not hold even over a DVR of characteristic $0$. 

\begin{ex} \label{Ex_1}
Let $R = {\mathbb{Z}}_{(p)}$ where $p$ is a prime in $\mathbb{Z}$,
$K=Qt(R)= \mathbb{Q}$ and $k=R/pR=\mathbb{Z}/p\mathbb{Z}$.
Let  $a=Y^p+Y+pX$ and $g=Z^p-a \in R[X,Y,Z]$. 
Then $R[X,Y,Z]=R[g]^{[2]}$; in particular,
$R[X,Y,Z]/(g) = R^{[2]}$. But $R[X,Y] \ne R[a]^{[1]}$.
\end{ex}

\begin{proof}
Let $Z'=Z-Y$. Then $R[X,Y,Z]=R[X,Y,Z']$ and $g={Z'}^p -pf(Z',Y)-Y-pX$ 
for some $f \in R[Z',Y]$. Let $D=R[g,Z']$.
We have $K[X,Y,Z]=K[g,Y,Z]=K[g,Z']^{[1]}$ and 
$k[X,Y,Z]=k[\bar{g},X,Z']=k[\bar{g},Z']^{[1]}$
where $\bar{g}$ denotes the image of $g$ in $k[X,Y,Z]$.
Since $p$ is prime in $R$, $p$ is prime in both $R[X,Y,Z]$ and $D$. 
Hence, by Theorem \ref{RS_Poly}, $R[X,Y,Z]=D^{[1]}=R[g]^{[2]}$. 
Let $\bar{a}$ denote the image of $a$ in $k[X,Y]$.
Since $k[\bar{a}]=k[Y+Y^p]$ is not algebraically closed in $k[X,Y]$, 
$\bar{a}$ is not a variable in $k[X,Y]$ and hence $a$ is not a variable in 
$R[X,Y]$.
\end{proof}

However the next result shows that Theorem \ref{Our_Th4} holds over any 
DVR $(R,t)$ of characteristic $0$ (without assuming
that the characterisic of $R/tR$ does not divide $n$),
if the element $a$ is such that 
$(R/tR)[\bar{a}]$ is algebraically closed in $(R/tR)[X,Y]$.

\begin{prop} \label{Prop_DVR_Alg_Closed}
Let $(R, t)$ be a DVR of characteristic $0$ with residue field $k$
and $g \in R[X,Y,Z](= R^{[3]})$ 
be of the form $g = bZ^n -a$ where $a,b \in R[X,Y]$, $b \ne 0$ and $n$ 
is an integer $\ge 2$. Suppose that $R[X,Y,Z]/(g) = R^{[2]}$ and 
$k[\bar{a}]$ is algebraically closed in $k[X,Y]$. Then 
$R[X,Y]=R[a]^{[1]}$ and $R[X,Y,Z] = R[Z,g]^{[1]}$.
\end{prop}

\begin{proof}
 
We see that $R[1/t][X,Y] = R[1/t][a]^{[1]}$ by Theorem \ref{Our_Th0a},
$t$ is prime in both $R[a]$ and $R[X,Y]$,
$tR[X,Y] \cap R[a] = tR[a]$ by Lemma \ref{Lem_Inert_R[a]} and
$(R/tR)[\bar{a}]$ is algebraically closed in $(R/tR)[X,Y]$ 
by hypothesis. Hence, by Theorem \ref{RS_Poly},  $R[X,Y]=R[a]^{[1]}$.
Let $B := R[X,Y,Z]/(g) (=R^{[2]})$ and denote the image 
of $Z$ in $B$ by $\z$. Then $B/(\z) = R[X,Y,Z]/(Z,  bZ^n-a) 
= R[X,Y]/(a) = R^{[1]}$ and hence, by the generalized epimorphism 
theorem of Bhatwadekar (\cite{B_GENEPI}, Theorem 3.7),
we have $B = R[\z]^{[1]}$. Let $C=R[Z]$. 
Identifying the image of $Z$ in $B$ with $Z$ itself, we 
have $C[X,Y]/(g)=C^{[1]}$. Since $C$ is a normal domain 
of characteristic $0$, again by Bhatwadekar's result 
(\cite{B_GENEPI}, Theorem 3.7), we have $C[X,Y] = C[g]^{[1]}$, 
i.e., $R[X,Y,Z] = R [g,Z]^{[1]}$.
\end{proof}

In view of Example \ref{Ex_1} and Proposition \ref{Prop_DVR_Alg_Closed}
we ask:

\begin{prob}
 Let $(R, t)$ be a DVR of characteristic 0 such that the 
characteristic of the residue field is positive, say $p$. 
Let $g = b Z^{pm} - a \in R[X,Y,Z]$ be such that 
 $R[X,Y,Z]/(g) = R^{[2]}$ where $a, b (\ne 0) \in R[X,Y]$ and $m \ge 1$. 
 Is then $R[X,Y,Z] = R[g]^{[2]}$? 
\end{prob}

\section{Planes of the form $bZ^n-a$ over rings containing a field} 
\label{sec_misc}

In this section we prove a generalized version of Theorem B 
(Theorem \ref{Our_Th2d}). 
The authors thank Neena Gupta for her observations on 
the earlier drafts of this paper which have resulted in the
formulation of Theorem \ref{Our_Th2d} in its present generality. 
We shall essentially follow
the approach of Bhatwadekar in (\cite{B_GENEPI}) and then
apply the result on residual variables (Theorem \ref{Res_Result}). 
We first state a result which will be needed in the proof of 
Theorem \ref{Our_Th2d}.
 
\begin{lem} \label{Lem_Val}
 Let $R$ be a Noetherian domain and let $b$ ($\ne 0$)
$\in R$. Then, for each non-zero prime ideal $P$ of $R$, there 
exists a discrete valuation ring $V$ with maximal ideal ${\m}_V$ 
together with a homomorphism 
$\phi: R \longrightarrow V$ such that $\phi(b) \ne 0$, 
$\phi^{-1}({\m}_V) = P$ and $V/{\m}_V$ is algebraic over $k(P)$.
\end{lem}
\begin{proof}
 Let $P$ be a non-zero prime ideal of $R$ and let $n$ be the height of $P$. 
 Since $R$ is a Noetherian domain, there exists a prime ideal $Q$ of $R$ 
of height $n-1$ such that $Q \subsetneqq P$ and $b \notin Q$. 
Let $D=R/Q$ and $\p =P/Q$ be the image of $P$ in $D$. 
Let $C$ be the normalisation of $D$ and $\mathcal{P}$ a prime ideal
of $C$ lying over $\p$. Set $V:=C_{\mathcal{P}}$ and 
${\m}_V:=\mathcal{P}V$, the maximal ideal of the local ring $V$.
Since the height of $\p$ (and hence that of $\mathcal{P}$) is one,
$V$ is a DVR. Now let $\phi$ denote the composite map $R \longrightarrow 
\ D(= R/Q) \longrightarrow C 
\longrightarrow V (=C_{\mathcal{P}})$.
Clearly, $\phi^{-1}({\m}_V) = P$, $\phi(b) \ne 0$ 
and $V/{\m}_V$ is algebraic over $k(P)$. 
\end{proof}

We now prove the main result of this section.

\begin{thm} \label{Our_Th2d}
Let $R$ be a Noetherian domain containing a field of 
characteristic $p \ge 0$ and $g \in R[X,Y,Z](=R^{[3]})$ 
be of the form $bZ^n-a$ where $a,b \in R[X, Y]$, $b \ne 0$ 
and $n$ is an integer $ \ge 2$ such that $p \nmid n$. Suppose 
that $R[X,Y,Z]/(g) = R^{[2]}$. Then 
$R[X, Y, Z] \otimes_R k(P) = (R[g, Z] \otimes_R k(P))^{[1]}$ 
and $R[X, Y] \otimes_R k(P) = (R[a] \otimes_R k(P))^{[1]}$ for all 
$P \in Spec(R)$. Thus, if $R$ contains $\mathbb Q$ or if $R$ is seminormal, 
then $R[X, Y, Z] = R[g, Z]^{[1]}$ and $R[X, Y] = R[a]^{[1]}$.
\end{thm}
\begin{proof}
Fix $P \in Spec(R)$. 
Let the images of $b$, $a$ and $g$ in 
$R[X,Y,Z] \otimes_R k(P)$ be $\bar{b}$, $\bar{a}$ and $\bar{g}$ respectively.
Let $k$ denote $k(P)$. We show that
$k[X,Y]=k[\bar{a}]^{[1]}$ and $k[X,Y,Z]=k[\bar{g}, Z]^{[1]}$.

If $ht \ P = 0$, we are done by Theorem \ref{Our_Th0a}. 
So we assume that
 $ht \ P = n \ge 1$. If $\bar{b} \ne 0$, then by Theorem \ref{Our_Th0a}, 
we are through. So we assume that $\bar{b} = 0$ (and hence $\bar{g} = \bar{a}$).

Using Lemma \ref{Lem_Val}, we have a DVR $(V,\pi)$ with a homomorphism 
$\phi: R \longrightarrow V$ such that $\phi(b) \ne 0$, 
$\phi^{-1}(\pi) = P$ and $V/(\pi)$ is algebraic over $k(P)$. 
Note that $V[X,Y,Z]/(g) = V^{[2]}$ and hence, by Theorem \ref{Our_Th4}, 
we have $V[X,Y,Z] = V[g,Z]^{[1]}$ and $V[X,Y] = V[a]^{[1]}$; in particular, 
$\frac{V}{(\pi)}[X,Y,Z] = \frac{V}{(\pi)} [\bar{g}, Z]^{[1]}$ and 
$\frac{V}{(\pi)}[X,Y] = \frac{V}{(\pi)}[\bar{a}]^{[1]}$. Now, since 
$\frac{k[X,Y]}{(\bar{a})}[Z] = k^{[2]}$, by Theorem \ref{AEH}, we have
$\frac{k[X,Y]}{(\bar{a})} = k^{[1]}$. Since $V/(\pi)$ is algebraic 
over $k$ and since $\frac{V}{(\pi)}[X,Y] = \frac{V}{(\pi)}[\bar{a}]^{[1]}$, 
by (\cite{Ganong_Plane_Curve}, Proposition 1.16), we have 
$k[X,Y]=k[\bar{a}]^{[1]}$ and hence $k[X,Y,Z] = k[\bar{g},Z]^{[1]}$.

Thus $R[X, Y, Z] \otimes_R k(P) = (R[g, Z] \otimes_R k(P))^{[1]}$ 
and $R[X, Y] \otimes_R k(P) = (R[a] \otimes_R k(P))^{[1]}$ for all 
$P \in Spec(R)$.

Now, if $R$ is seminormal or contains $\mathbb Q$,
then $R[X, Y, Z] = R[g, Z]^{[1]}$ and $R[X, Y] = R[a]^{[1]}$
by Theorem \ref{Res_Result}.
\end{proof}

\begin{rem}
(1) If $R$ is seminormal or $R$ contains $\mathbb Q$,
then under the hypotheses of 
Theorem \ref{Our_Th2d}, one can show, by suitable 
reductions, that $R[X, Y, Z] = R[g, Z]^{[1]}$ and 
$R[X, Y] = R[a]^{[1]}$, even when $R$ is non-Noetherian.

(2) If $R$ is a UFD, then under the hypotheses of 
Theorem \ref{Our_Th2d}, $b \in R[X_0]$ where $X_0 \in R[X,Y]$
and $K[X,Y]=K[X_0, a]$. This follows from the proof of Theorem 
\ref{Our_Th4}.
\end{rem}

\noindent
{\bf Acknowledgement:} The authors thank S. M. Bhatwadekar, 
Neena Gupta and N. Onoda for many fruitful suggestions.

\bibliography{reference}
\bibliographystyle{amsalpha}

\end{document}